\documentclass{amsart}
\usepackage{a4}

\usepackage[dvips]{epsfig}
\usepackage{latexsym}
\usepackage{amssymb}
\usepackage{algorithm}
\usepackage{algorithmic}
\usepackage{amsfonts}
\usepackage{amsmath}
\usepackage{indentfirst}

\title{Permutations defining convex permutominoes}

\newtheorem{proposition}{Proposition}
\newtheorem{corollary}{Corollary}
\newtheorem{example}{Example}
\newtheorem{theorem}{Theorem}
\author{A. Bernini
 \and F. Disanto
 \and
 R. Pinzani
 \and
 S. Rinaldi
}

\address{{\rm Antonio Bernini, Renzo Pinzani:} \newline Universit\`a di Firenze, Dipartimento di
Sistemi e Informatica \newline viale Morgagni 65, 50134 Firenze, Italy {\tt (bernini, pinzani@dsi.unifi.it}).}

\address{{\rm Filippo Disanto, Simone Rinaldi:} \newline Universit\`a di Siena,
Dipartimento di Scienze Matematiche e Informatiche \newline Pian dei Mantellini 44, 53100 Siena, Italy ({\tt
rinaldi@unisi.it}).}

\date{} 

\begin{document}\date{}
\maketitle
\begin{abstract}
A permutomino of size $n$ is a polyomino determined by particular pairs $(\pi_1, \pi_2)$ of permutations of size
$n$, such that $\pi_1(i)\neq \pi_2(i)$, for $1\leq i\leq n$. Here we determine the combinatorial properties and,
in particular, the characterization for the permutations defining convex permutominoes.

Using such a characterization, these permutations can be uniquely represented in terms of the so called square
permutations, introduced by Mansour and Severini. Then, we provide a closed formula for the number of these
permutations with size $n$.

\end{abstract}

\section{Convex polyominoes}

\begin{figure}[htb]
\begin{center}
\epsfig{file=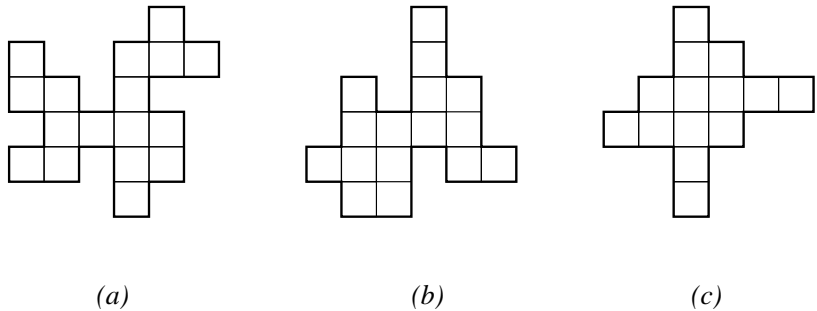}
\end{center}
\caption{(a) a polyomino; (b) a column convex polyomino which is not row convex; (c) a convex polyomino.
\label{polyom}}
\end{figure}

In the plane $\Bbb Z \times \Bbb Z$ a {\em cell} is a unit square, and a {\em polyomino} is a finite connected
union of cells having no cut point. Polyominoes are defined up to translations (see Figure~\ref{polyom}~(a)). A
{\em column} ({\em row}) of a polyomino is the intersection between the polyomino and an infinite strip of cells
lying on a vertical (horizontal) line.

Polyominoes were introduced by Golomb~\cite{dbintr39}, and then they have been studied in several mathematical
problems, such as tilings~\cite{Gi, Go}, or games~\cite{Ga} among many others. The enumeration problem for
general polyominoes is difficult to solve and still open. The number $a_n$ of polyominoes with $n$ cells is known
up to $n=56$~\cite{jensen-guttmann} and asymptotically, these numbers satisfy the relation $\smash{\lim_n
\left(a_n \right)^{1/n}=\mu}$, \ $3.96 < \mu <4.64$, where the lower bound is a recent improvement
of~\cite{bmrr}.

\medskip

In order to simplify enumeration problems of polyominoes, several subclasses were defined by combining the two
simple notions of {\em convexity} and {\em directed growth}. A polyomino is said to be {\em column convex} (resp.
{\em row convex}) if every its column (resp. row) is connected (see Figure~\ref{polyom}~$(b)$). A polyomino is
said to be {\em convex}, if it is both row and column convex (see Figure~\ref{polyom}~$(c)$). The {\em area} of a
polyomino is just the number of cells it contains, while its {\em semi-perimeter} is half the number of edges of
cells in its boundary. Thus, for any convex polyomino the semi-perimeter is the sum of the numbers of its rows
and columns. Moreover, any convex polyomino is contained in a rectangle in the square lattice which has the same
semi-perimeter, called {\em minimal bounding rectangle}.

\medskip

A significant result in the enumeration of convex polyominoes was first obtained by Delest and Viennot in
\cite{DV}, where the authors proved that the number $\ell_n$ of convex polyominoes with semi-perimeter equal to
$n+2$ is:
\begin{equation}\label{eq1}
\ell_{n+2}=(2n+11)4^n \, - \, 4(2n+1) \, {2n \choose n}, \quad n \geq 2; \quad \ell_{0}=1, \quad \ell_{1}=2.
\end{equation}
\noindent This is sequence $A005436$ in \cite{sloane}, the first few terms being:
$$1,2,7,28,120,528,2344,10416, \ldots .$$

During the last two decades convex polyominoes, and several combinatorial objects obtained as a generalizations
of this class, have been studied by various points of view. For the main results concerning the enumeration and
other combinatorial properties of convex polyominoes we refer to \cite{mbm2, mbm, gutman, chang}.

\medskip

There are two other classes of convex polyominoes which will be useful in the paper, the {\em directed convex
polyominoes} and the {\em parallelogram}. A polyomino is said to be {\em directed} when each of its cells can be
reached from a distinguished cell, called the root, by a path which is contained in the polyomino and uses only
north and east unitary steps.

\begin{figure}[htb]
\centerline{\hbox{\psfig{figure=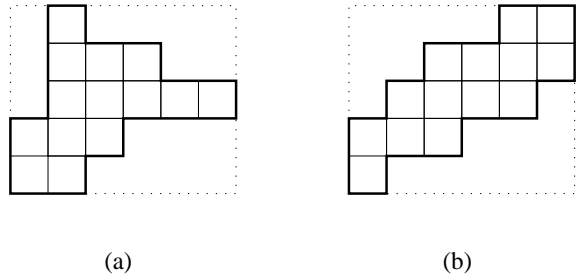}}} \caption{(a)  A directed convex polyomino; (b) a
parallelogram polyomino.\label{dirconv}}
\end{figure}

A polyomino is {\em directed convex} if it is both directed and convex (see Figure~\ref{dirconv}~(a)). It is
known that the number of directed convex polyominoes of semi-perimeter $n+2$ is equal to the $n$th central
binomial coefficient, i.e.,
\begin{equation}\label{bin}
b_n={2n \choose n},
\end{equation}
sequence A000984 in \cite{sloane}.

Finally, {\em parallelogram} polyominoes are a special subset of the directed convex ones, defined by two lattice
paths that use north and east unit steps, and intersect only at their origin and extremity. These paths are
called the {\em upper} and the {\em lower path} (see Figure~\ref{dirconv}~(b)). It is known~\cite{stan} that the
number of parallelogram polyominoes having semi-perimeter $n+1$ is the $n$-th {\em Catalan number} (sequence
A000108 in \cite{sloane}),
\begin{equation}\label{cat}
c_n=\frac{1}{n+1} {2n \choose n}.
\end{equation}


\section{Convex permutominoes}

\begin{figure}[htb]
\begin{center}
\epsfig{file=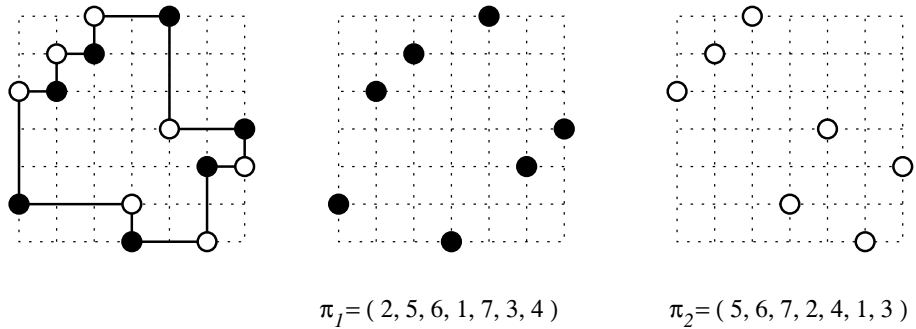} \caption{A permutomino and the two associated permutations. \label{permu1}}
\end{center}
\end{figure}

Let $P$ be a polyomino without holes, having $n$ rows and columns, $n\geq 1$; we assume without loss of
generality that the south-west corner of its minimal bounding rectangle is placed in $(1,1)$. Let ${\mathcal A}=
\left( A_1, \ldots ,A_{2(r+1)} \right)$ be the list of its vertices (i.e., corners of its boundary) ordered in a
clockwise sense starting from the lowest leftmost vertex. We say that $P$ is a {\em permutomino} if ${\mathcal
P}_1 = \left( A_1, A_3, \ldots ,A_{2r+1} \right)$ and ${\mathcal P}_2 = \left( A_2, A_4, \ldots ,A_{2r+2}
\right)$ represent two permutations of ${\mathcal S}_{n+1}$, where, as usual, ${\mathcal S}_{n}$ is the symmetric
group of size $n$. Obviously, if $P$ is a permutomino, then $r=n$, and $n+1$ is called the {\em size} of the
permutomino. The two permutations defined by ${\mathcal P}_1$ and ${\mathcal P}_2$ are indicated by $\pi _1(P)$
and $\pi _2(P)$, respectively (see Figure \ref{permu1}).

From the definition any permutomino $P$ has the property that, for each abscissa (ordinate) there is exactly one
vertical (horizontal) side in the boundary of $P$ with that coordinate. It is simple to observe that this
property is also a sufficient condition for a polyomino to be a permutomino. By convention we also consider the
empty permutomino of size $1$, associated with $\pi=(1)$.

Permutominoes were introduced by F.~Incitti in \cite{incitti} while studying the problem of determining the
$\widetilde{R}$-polynomials (related with the Kazhdan-Lusztig R-polynomials) associated with a pair
$(\pi_1,\pi_2)$ of permutations. Concerning the class of polyominoes without holes, our definition (though
different) turns out to be equivalent to Incitti's one, which is more general but uses some algebraic notions not
necessary in this paper.

\medskip

Let us recall the main enumerative results concerning convex permutominoes. In \cite{fanti}, using bijective
techniques, it was proved that the number of {\em parallelogram permutominoes} of size $n+1$ is equal to $c_n$
and that the number of {\em directed-convex permutominoes} of size $n+1$ is equal to $ \frac{1}{2} \, b_n $,
where, throughout all the paper, $c_n$ and $b_n$ will denote, respectively, the Catalan numbers and the central
binomial coefficients. Finally, in \cite{rinaldi} it was proved, using the ECO method, that the number of {\em
convex permutominoes} of size $n+1$ is:
\begin{equation}\label{co}
 \, 2 \, (n+3) \, 4^{n-2} \, - \, \frac{n}{2} \, {{2n} \choose {n}} \qquad n\geq 1.
\end{equation}
The first terms of the sequence are
$$ 1,1,4,18,84,394,1836,8468, \ldots $$ (sequence A126020) in \cite{sloane}).
The same formula has been obtained independently by Boldi et al. in \cite{milanesi}. The main results concerning
the enumeration of classes of convex permutominoes are listed in Table~\ref{tttt}, where the first terms of the
sequences are given starting from $n=1$, and are taken from~\cite{rinaldi,fanti}.

\begin{table}\label{tttt}
\begin{tabular}{l|l|l}
  $Class$ & {\em First terms } & {\em Closed form/rec. relation} \\
      & & \\
  \hline
      & &  \\
  convex &$1,1,4,18,84,394, \ldots $  &$C_{n+1} = 2 \, (n+3) \, 4^{n-2} \, - \,
                                                        \frac{n}{2} \, {{2n} \choose {n}}$ \\
      & & \\
  \hline
      & &  \\
 \hspace{-.35cm} \begin{tabular}{l}
  directed \\
  convex
  \end{tabular}
  &$1,1,3,10,35,126, \ldots $ &$D_{n+1} = \frac{1}{2} \, b_{n}$ \\
      & & \\
  \hline
      & &  \\
  parallelogram &$1,1,2,5,14,42,132, \ldots $ &$P_{n+1} =  c_{n}$ \\
  & & \\
  \hline
  & &  \\
 \hspace{-.35cm} \begin{tabular}{l}
  symmetric \\
  (w.r.t. $x=y$)
  \end{tabular}
  &$1,1,2,4,10,22,54, \ldots $
  &\hspace{-.3cm} $\begin{array}{l} S_{n+1} =  (n+3)2^{n-2}-n{{n-1}\choose{\lfloor{\frac{n-1}{2}}\rfloor}} \\
  \\
  \phantom{S_{n+1} = \,} -(n-1){{n-2}\choose{\lfloor{\frac{n-2}{2}}\rfloor}}
  \end{array}$\\
  & & \\
  \hline

 & &  \\
  centered &$1,1,4,16,64,256, \ldots $   &$Q_{n} = 4^{n-2}$ \\
 & &  \\
\hline
 & &  \\
  bi-centered &$1,1,4,14,48,164, \ldots $   &$T_n \, = \, 4 T_{n-1} \, - \, 2 T_{n-2},
                                            \quad n\geq 3$ \\
 & &  \\
\hline
 & &  \\
  stacks &$1,1,2,4,8,16,32, \ldots $   &$W_n \, = \, 2^{n-2}$ \\
 & &  \\
\end{tabular}
\end{table}

\bigskip

\bigskip

\noindent {\bf Notation.} Throughout the whole paper we are going to use the following notations:\begin{itemize}
    \item $\mathcal C_n$ is the set of convex permutominoes of size $n$;
    \item $C_n$ is the cardinality of $\mathcal C_n$;
    \item $C(x)$ is the generating function of the sequence $\{C_n\}_{n\geq 2}$.
\end{itemize}
Moreover, if $\pi$ is a permutation of size $n$, then we define its \emph{reversal} $\pi^R$ and its
\emph{complement} $\pi^C$ as follows: $\pi^R(i)=\pi(n+1-i)$ and $\pi^C(i)=n+1-\pi(i)$, for each $i=1,\ldots,n$.


\section{Permutations associated with convex permutominoes}

\begin{figure}[htb]
\begin{center}
\epsfig{file=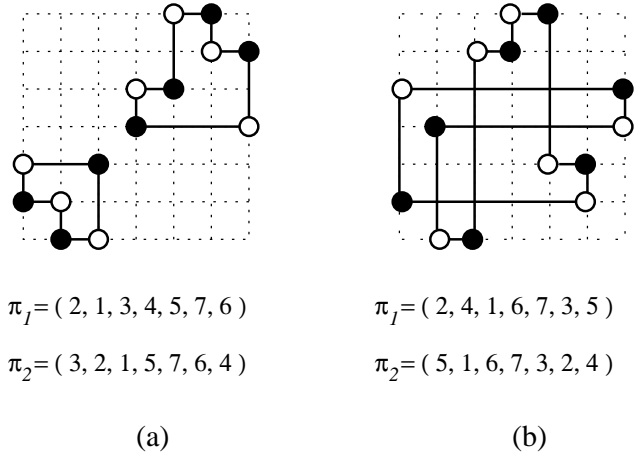} \caption{Two permutations $\pi _1$ and $\pi _2$ of ${\mathcal S}_{n}$, satisfying 1 and
2, do not necessarily define a permutomino, since two problems may occur: (a) two disconnected sets of cells; (b)
the boundary crosses itself. \label{permuz}}
\end{center}
\end{figure}

Given a permutomino $P$, the two permutations we associate with $P$ are denoted by $\pi _1$ and $\pi _2$ (see
Figure \ref{permu1}). While it is clear that any permutomino of size $n\geq 2$ uniquely determines two
permutations $\pi _1$ and $\pi _2$ of ${\mathcal S}_{n}$, with
\begin{description}
    \item[1] $\pi_1 (i) \neq \pi_2 (i)$, $1 \leq i \leq n$,
    \item[2] $\pi_1(1)<\pi_2(1)$, and $\pi_1(n)>\pi_2(n)$,
\end{description}
not all the pairs of permutations $\left ( \pi _1 , \pi _2 \right )$ of $n$ satisfying 1  and 2 define a
permutomino: Figure \ref{permuz} depicts the two problems which may occur.

In \cite{fanti} the authors give a simple constructive proof that every permutation of ${\mathcal S}_n$ is
associated with at least one column convex permutomino.

\begin{proposition}\label{ovo}
If $\pi \in {\mathcal S}_n$ then there is at least one column convex permutomino $P$ such that $\pi = \pi_1(P)$
or $\pi = \pi_2 (P)$.
\end{proposition}

\begin{figure}[htb]
\begin{center}
\epsfig{file=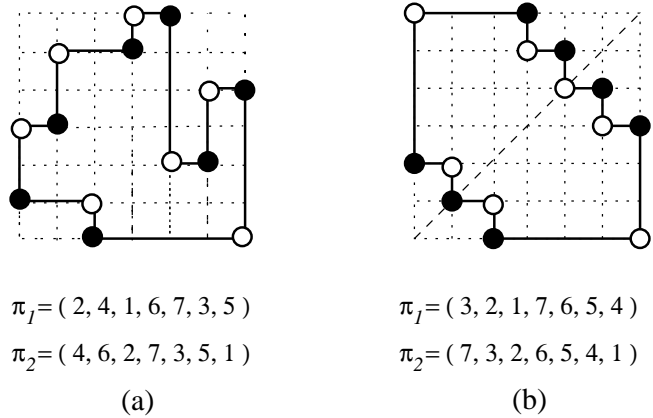} \caption{(a) a column convex permutomino associated with the permutation $\pi_1$ in
Figure \ref{permuz} (b); (b) the symmetric permutomino associated with the involution $\pi_1=(3,2,1,7,6,5,4)$.
\label{casesb}}
\end{center}
\end{figure}

For instance, Figure \ref{casesb} (a) depicts a column convex permutomino associated with the permutation $\pi_1$
in Figure \ref{permuz} (b).

\begin{figure}[htb]
\begin{center}
\epsfig{file=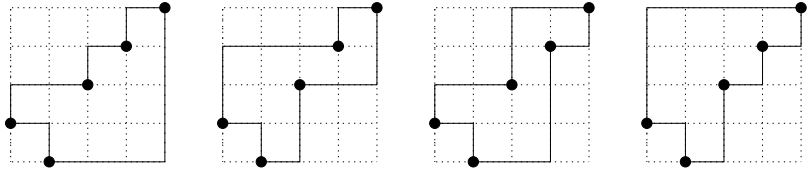,width=4.2in,clip=} \caption{The four convex permutominoes associated with $(2,1,3,4,5)$.
\label{pippe}}
\end{center}
\end{figure}

The statement of Proposition~\ref{ovo} does not hold for convex permutominoes. Therefore, in this paper we
consider the class ${\mathcal C}_n$ of convex permutominoes of size $n$, and study the problem of giving a
characterization for the set of permutations defining convex permutominoes,
$$ \left \{ \, (\pi_1(P),\pi_2(P)) \, : \, P \in {\mathcal C}_n \, \right \}.$$
 Moreover, let us consider the following
subsets of $S_n$:
$$ \widetilde{\mathcal C}_n=\{ \, \pi_1 (P) \, : \, P \in {\mathcal C}_n \,
    \}, \qquad  \widetilde{\mathcal C}'_n=\{ \, \pi_2 (P) \, : \, P \in {\mathcal C}_n \,
    \}. $$
It is easy to prove the following properties:
\begin{enumerate}
    \item $\left | \widetilde{\mathcal C}_n \right | = \left |  \widetilde{\mathcal C}'_n \right
    |$,
    \item $\pi \in \widetilde{\mathcal C}_n$ if and only if $\pi^R \in \widetilde{\mathcal
    C}'_n$.
    \item If $P$ is symmetric according to the diagonal $x=y$, then $\pi _1(P)$ and
    $\pi _2(P)$ are both {\em involutions} of ${\mathcal S}_n$. We recall that an involution
    is a permutation where all the cycles have length at most $2$ (see for instance Figure \ref{casesb} (b)).
    Figures \ref{pippe} and \ref{classi} show permutominoes where only $\pi_1$ is an involution, and this
    condition is not sufficient for the permutomino to be symmetric.
\end{enumerate}

Given a permutation $\pi \in {\mathcal S}_n$, we say that $\pi$ is {\em $\pi_1$-associated} (briefly {\em
associated}) with a permutomino $P$, if $\pi = \pi _1(P)$. With no loss of generality, we will study the
combinatorial properties of the permutations of $\widetilde{\mathcal C}_n$, and we will give a simple way to
recognize if a permutation $\pi$ is in $\widetilde{\mathcal C}_n$ or not. Moreover, we will study the cardinality
of this set. In particular, we will exploit the relations between the cardinalities of ${\mathcal C}_n$ and of
$\widetilde{\mathcal C}_n$.

\medskip
\noindent For small values of $n$ we have that:
\begin{eqnarray*}
  \widetilde{\mathcal C}_1 &=& \{ 1 \}, \\
  \widetilde{\mathcal C}_2 &=& \{ 12 \}, \\
  \widetilde{\mathcal C}_3 &=& \{ 123, 132, 213 \}, \\
  \widetilde{\mathcal C}_4 &=& \{ 1234, 1243, 1324 , 1342, 1423,1432, 2143,   \\
   & &  \, \, \, 2314, 2134, 2413, 3124, 3142, 3214  \}.\\
\end{eqnarray*}
As a main result we will prove that the cardinality of $\widetilde{\mathcal C}_{n+1}$ is
\begin{equation}\label{cnv}
 \, 2 \, (n+2) \, 4^{n-2} \, - \, \frac{n}{4} \, \left ( {\frac{3-4n}{1-2n}}
\right ) \, {{2n} \choose {n}}, \qquad n \geq 1.
\end{equation}
defining the sequence $1,1,3,13,62,301,1450, \ldots$, recently added to \cite{sloane} as A122122. For any $\pi
\in \widetilde{\mathcal C}_n$, let us consider also
$$[\pi]= \{ P\in {\mathcal C}_n : \pi _1 (P)=\pi \},$$
i.e., the set of convex permutominoes associated with $\pi$. For instance, there are $4$ convex permutominoes
associated with $\pi = (2,1,3,4,5)$, as depicted in Figure~\ref{pippe}. In this paper we will also give a simple
way of computing $[\pi ]$, for any given $\pi \in \widetilde{\mathcal C}_n$.

\subsection{A matrix representation of convex permutominoes}\label{rif}

Before going on with the study of convex permutominoes, we would like to point out a simple property of their
boundary, related to {\em reentrant} and {\em salient points}. Let us briefly recall the definition of these
objects.

\begin{figure}[htb]
\begin{center}
\epsfig{file=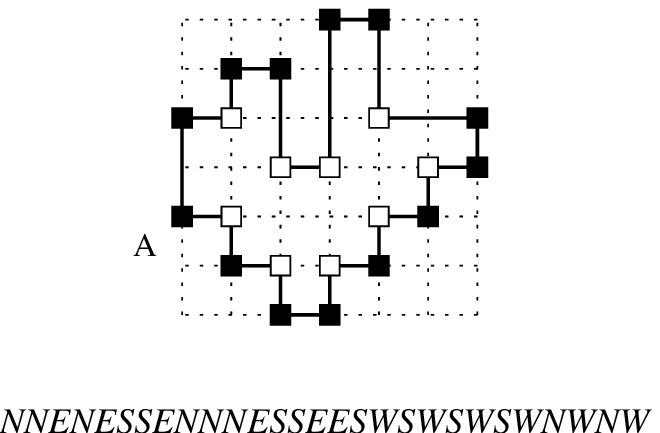} \caption{The coding of the boundary of a polyomino, starting from $A$ and moving in a
clockwise sense; its salient (resp. reentrant) points are indicated by black (resp. white) squares.
\label{reentrant}}
\end{center}
\end{figure}

Let $P$ be a polyomino; starting from the leftmost point having minimal ordinate, and moving in a clockwise
sense, the boundary of $P$ can be encoded as a word in a four letter alphabet, $\{ N, E, S, W \}$, where $N$
(resp., $E$, $S$, $W$) represents a {\em north} (resp., {\em east}, {\em south}, {\em west}) unit step. Any
occurrence of a sequence $NE$, $ES$, $SW$, or $WN$ in the word encoding $P$ defines a {\em salient point} of $P$,
while any occurrence of a sequence $EN$, $SE$, $WS$, or $NW$ defines a {\em reentrant point} of $P$ (see for
instance, Figure \ref{reentrant}).

In \cite{daurat} and successively in \cite{brlek}, in a more general context, it was proved that in any polyomino
the difference between the number of salient and reentrant points is equal to $4$.

In a convex permutomino of size $n+1$ the length of the word coding the boundary is $4n$, and we have $n+3$
salient points and $n-1$ reentrant points; moreover we observe that a reentrant point cannot lie on the minimal
bounding rectangle. This leads to the following remarkable property:

\begin{proposition}\label{spigoli}
The set of reentrant points of a convex permutomino of size $n+1$ defines a permutation matrix of dimension
$n-1$, $n\geq 1$.
\end{proposition}

For simplicity of notation, we agree to group the reentrant points of a convex permutomino in four classes; in
practice we choose to represent the reentrant point determined by a sequence $EN$ (resp. $SE$, $WS$, $NW$) with
the symbol $\alpha$ (resp. $\beta$, $\gamma$, $\delta$).
\begin{figure}[htb]
\begin{center}
\epsfig{file=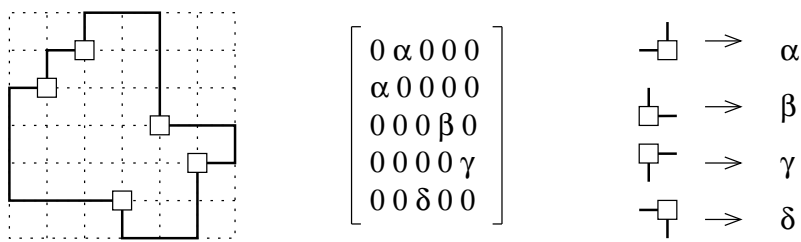} \caption{The reentrant points of a convex permutomino uniquely define a permutation
matrix in the symbols $\alpha$, $\beta$, $\gamma$ and $\delta$. \label{permu2}}
\end{center}
\end{figure}

Using this notation we can state the following simple characterization for convex permutominoes:
\begin{proposition}\label{car_mat}
A convex permutomino of size $n\geq 2$ is uniquely represented by the permutation matrix defined by its reentrant
points, which has dimension $n-2$, and uses the symbols $\alpha , \, \beta , \, \gamma , \, \delta $, and such
that for all points $A,B,C,D$, of type $\alpha$, $\beta$, $\gamma$ and $\delta$, respectively, we have:
\begin{enumerate}
    \item $x_A < x_B$, $x_D < x_C$, $y_A > y_D$, $y_B > y_C$;
    \item $\neg (x_A > x_C \, \wedge \, y_A < y_C )$ and
    $\neg (x_B < x_D \, \wedge \, y_B < y_D )$,
    \item the ordinates of the $\alpha$ and of $\gamma$ points are strictly increasing, from left to right;
    the ordinates of the $\beta$ and of $\delta$ points are strictly decreasing, from left to right.
\end{enumerate}
where $x$ and $y$ denote the abscissa and the ordinate of the considered point.
\end{proposition}
\begin{figure}[htb]
\begin{center}
\epsfig{file=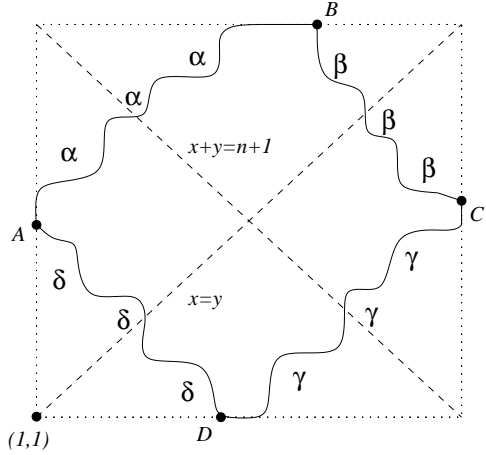,width=2.5in,clip=} \caption{A sketched representation of the $\alpha$, $\beta$, $\gamma$ and
$\delta$ paths in a convex permutomino. \label{d}}
\end{center}
\end{figure}
Just to give a more informal explanation, on a convex permutomino, let us consider the special points
$$A=(1,\pi_1(1)), \quad B=(\pi_1^{-1}(n),n), \quad C=(n,\pi_1(n)),
\quad D=(\pi_1^{-1}(1),1).$$ The path that goes from $A$ to $B$ (resp. from $B$ to $C$, from $C$ to $D$, and from
$D$ to $A$) in a clockwise sense is made only of $\alpha$ (resp. $\beta$, $\gamma$, $\delta$) points, thus it is
called the $\alpha$-path (resp. $\beta$-path, $\gamma$-path, $\delta$-path) of the permutomino. The situation is
schematically sketched in Figure~\ref{d}.

\medskip

\noindent From the characterization given in Proposition \ref{car_mat} we have the following two properties:

\begin{description}
    \item[(z1)] the $\alpha $ points are never below the diagonal $x=y$, and the $\gamma$ points
    are never above the diagonal $x=y$.
    \item[(z2)] the $\beta $ points are never below the diagonal $x+y=n+1$, and the $\delta$ points
    are never above the diagonal $x+y=n+1$.
\end{description}

\subsection{Characterization and combinatorial properties of $\widetilde{\mathcal C}_n$}

Let us consider the problem of establishing, for a given permutation $\pi \in {\mathcal S}_n$, if there is at
least a convex permutomino $P$ of size $n$ such that $\pi_1 (P)=\pi$.

\noindent Let $\pi$ be a permutation of ${\mathcal S}_n$, we define $\mu(\pi)$ (briefly $\mu$) as the maximal
upper unimodal sublist of $\pi$ ($\mu$ retains the indexing of $\pi$).

\medskip

\noindent Specifically, if $\mu$ is denoted by $\left ( \mu (i_1 ) , \ldots , n , \ldots , \mu (i_m) \right ),$
then we have the following:
\begin{enumerate}
    \item $\mu (i_1)=\mu(1)=\pi(1)$;
    \item if $n\notin \{ \mu (i_1), \ldots , \mu (i_k) \}$, then $\mu (i_{k+1})=\pi (i_{k+1})$ such that
    \begin{description}
        \item[i] $i_k < i < i_{k+1}$ implies $\pi(i)<\mu (i_k)$, and
        \item[ii] $\pi (i_{k+1}) > \mu (i_k)$;
    \end{description}
    \item if $n\in \{ \mu (i_1), \ldots , \mu (i_k) \}$, then $\mu (i_{k+1})=\pi (i_{k+1})$ such that
    \begin{description}
        \item[i] $i_k < i < i_{k+1}$ implies $\pi(i)<\pi (i_{k+1})$, and
        \item[ii] $\pi (i_{k+1}) < \mu (i_k)$.
    \end{description}
\end{enumerate}

\noindent Summarizing we have:
$$ \mu(i_1)=\mu(1)=\pi(1) < \mu (i_2) < \ldots < n > \ldots \mu (i_m)=\mu(n)=\pi(n). $$

\begin{figure}[htb]
\begin{center}
\centerline{\hbox{\includegraphics[width=0.75\textwidth]{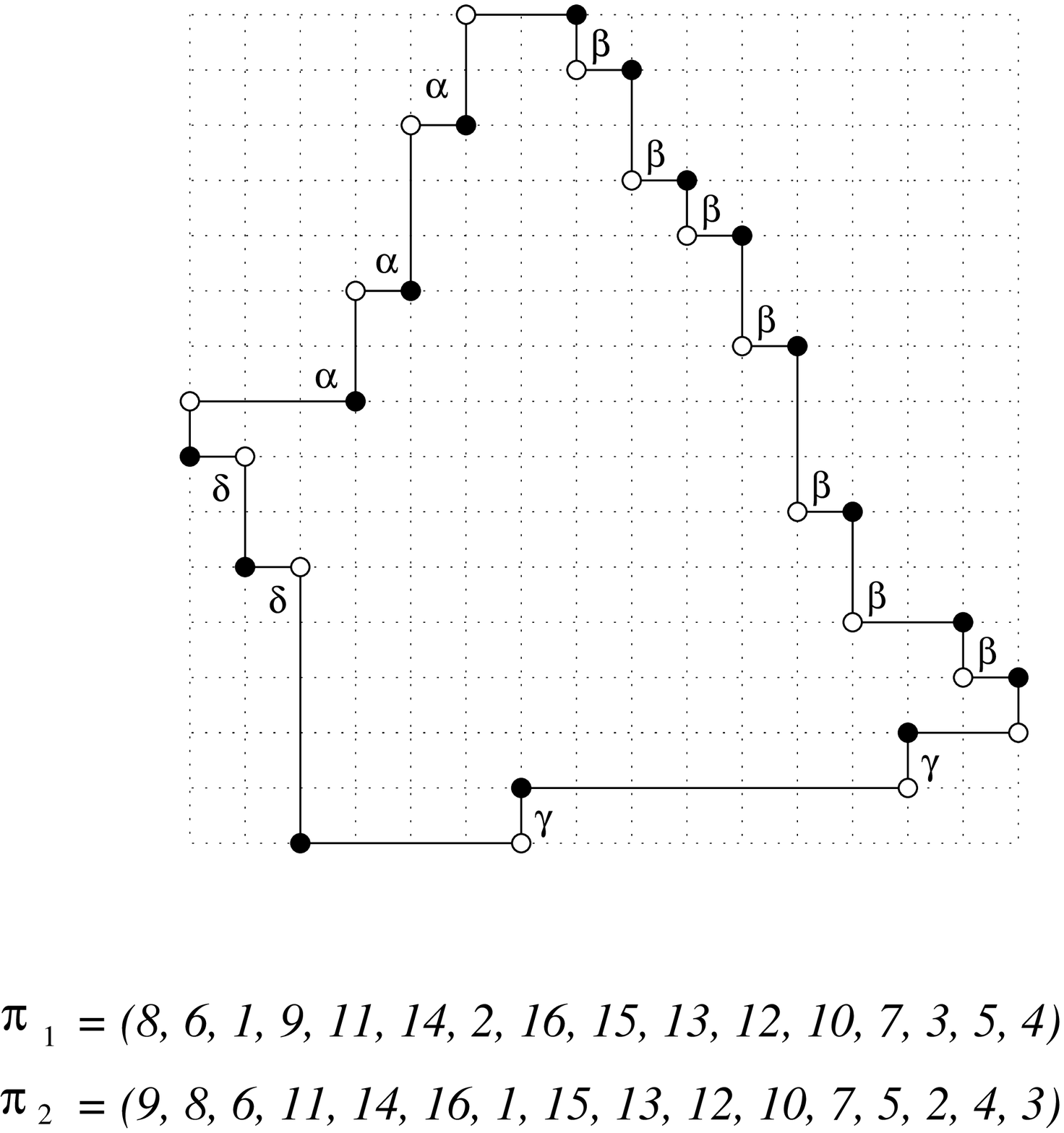}}} \caption{A convex permutomino and the
associated permutations.}\label{conv_new}
\end{center}
\end{figure}

\noindent Moreover, let $\sigma (\pi)$ (briefly $\sigma$) denote $\left(\sigma (j_1) , \ldots , \sigma (j_r)
\right )$ where:
\begin{enumerate}
    \item $\sigma (j_1)=\sigma(1)=\pi(1)$, $\sigma (j_r)=\sigma (n)=\pi (n)$, and
    \item if $1<j_k<j_r$, then $\sigma(j_k)=\pi(j_k)$ if and only if \\
    $\pi(j_k) \notin \{ \mu(i_1), \ldots , \mu (i_m) \}$.
\end{enumerate}

\medskip

We note that the sequence $\mu$ can be defined in terms of left-right and right left-maxima. A \emph{left-right
maximum} (resp. \emph{right-left maximum}) of a given permutation $\tau$ is an entry $\tau(j)$ such that
$\tau(j)>\tau(i)$ for each $i<j$ (for each $i>j$). Let $u=(u_{i_1},u_{i_2},\ldots ,u_{i_s})$ be the sequence of
the left-right maxima of $\pi$ with $u_{i_1}=\pi(1)<u_{i_2}<\ldots<u_{i_s}=n$, and let
$v=(v_{j_1},v_{j_2},\ldots,v_{j_t})$ be the sequence of the right-left maxima (read from the left) with
$v_{j_1}=n>v_{j_2}>\ldots>v_{j_t}=\pi(n)$. The sequence $\mu$ coincides with the sequence obtained by connecting
$u$ with $v$, observing that, clearly, $u_{i_s}=v_{j_1}=n$. In other words it is $\mu=(u_{i_1},u_{i_2},\ldots
,u_{i_s}(=v_{j_1}),v_{j_1},v_{j_2},\ldots,v_{j_t})$.

\begin{example}
{\em Consider the convex permutomino of size $16$ represented in Fig.~\ref{conv_new}. We have
$$ \pi_1 = (8,6,1,9,11,14,2,16,15,13,12,10,7,3,5,4), $$
\noindent and we can determine the decomposition of $\pi$ into the two subsequences $\mu$ and~$\sigma$:

\begin{center}
\begin{tabular}{c|cccccccccccccccc}
   &$1$ &$2$ &$3$ &$4$ &$5$ &$6$ &$7$ &$8$ &$9$ &$10$ &$11$ &$12$ &$13$ &$14$ &$15$ &$16$ \\
  \hline
  $\mu$ & $8$ & - & - & $9$ & $11$ &$14$ &- &$16$ &$15$ &$13$ &$12$ &$10$ &$7$ &- & $5$ & $4$ \\
  $\sigma$ &$8$ & $6$ & $1$ & - & - & - & 2 & - & - & - & - & - & - &$3$ & - &$4$ \\
\end{tabular}
\end{center}

For the sake of brevity, when there is no possibility of misunderstanding, we use to represent the two sequences
omitting the empty spaces, as
$$\mu = (8,9,11,14,16,15,13,12,10,7,5,4), \quad
\sigma=(8,6,1,2,3,4).$$}
\end{example}
While $\mu$ is upper unimodal by definition, here $\sigma$ turns out to be lower unimodal. In fact from the
characterization given in Proposition \ref{car_mat} we have that
\begin{proposition}\label{und}
If $\pi$ is associated with a convex permutomino then the sequence $\sigma$ is lower unimodal.
\end{proposition}

In this case, similarly to the sequence $\mu$, also the sequence $\sigma$ can be defined in terms of left-right
and right-left minima. A \emph{left-right minimum} (resp. \emph{right-left minimum}) of a given permutation
$\tau$ is an entry $\tau(j)$ such that $\tau(j)<\tau(i)$ for each $i<j$ (for each $i>j$). If $\sigma$ is lower
unimodal, then it is easily seen to be the sequence of the left-right minima followed by the sequence of the
right-left minima (read from the left), recalling that the entry $1$ is both a left-right minimum and a
right-left minimum.

\medskip

\begin{figure}[htb]
\begin{center}
\centerline{\hbox{\includegraphics[width=0.90\textwidth]{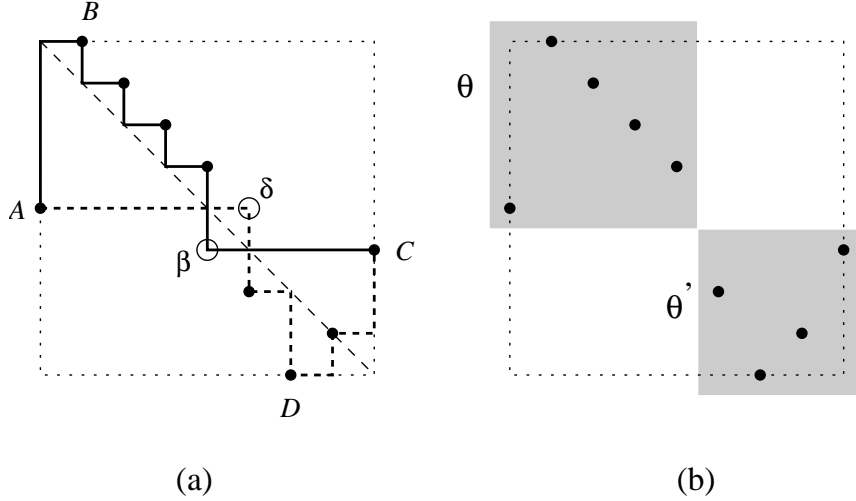}}} \caption{ (a) there is no convex
permutomino associated with $\pi =(5,9,8,7,6,3,1,2,4)$, since $\sigma$ is lower unimodal but the $\beta$ path
passes below the diagonal $x+y=10$. The $\beta$ point below the diagonal and the corresponding $\delta$ point
above the diagonal are encircled. (b) The permutation $\pi =(5,9,8,7,6,3,1,2,4)$ is the direct difference $\pi =
(1,5,4,3,2) \ominus (3,2,1,4)$.}\label{conv_news}
\end{center}
\end{figure}

The conclusion of Proposition \ref{und} is a necessary condition for a permutation $\pi$ to be associated with a
convex permutomino, but it is not sufficient. For instance, if we consider the permutation
$\pi=(5,9,8,7,6,3,1,2,4)$, then $\mu = (5,9,8,7,6,4)$, and $\sigma=(5,3,1,2,4)$ is lower unimodal, but as shown
in Figure \ref{conv_news} (a) there is no convex permutomino associated with $\pi$. In fact any convex
permutomino associated with such a permutation has a $\beta$ point below the diagonal $x+y=10$ and,
correspondingly, a $\delta$ point above this diagonal. Thus the $\beta$ and the $\delta$ paths cross themselves.

\medskip

In order to give a necessary and sufficient condition for a permutation $\pi$ to be in $\widetilde{\mathcal
C}_n$, let us recall that, given two permutations $\theta=(\theta_1, \ldots ,\theta_m) \in {\mathcal S}_m$ and
$\theta'=(\theta'_1, \ldots ,\theta'_{m'})\in {\mathcal S}_{m'}$, their {\em direct difference} $\theta \ominus
\theta '$ is a permutation of ${\mathcal S}_{m+m'}$ defined as
$$ (\theta_1+m', \ldots ,\theta_m+m',\theta'_1, \ldots ,\theta'_{m'} ). $$
A pictorial description is given in Figure~\ref{conv_news}~(b), where $\theta = (1,5,4,3,2)$, $\theta'=
(3,2,1,4)$, and their direct difference is $\theta \ominus \theta' =(5,9,8,7,6,3,1,2,4)$ .

\medskip

Finally the following characterization holds.
\begin{theorem}\label{caratt_conv}
Let $\pi \in {\mathcal S}_n$ be a permutation. Then $\pi \in \widetilde{\mathcal C}_n$ if and only if:

\begin{enumerate}
    \item $\sigma$ is lower unimodal, and
    \item there are no two permutations, $\theta \in\theta_m $ , and $\theta' \in \theta_m'$, such that
    $m+m'=n$, and $\pi = \theta \ominus \theta'$.
\end{enumerate}
\end{theorem}

\noindent {\bf (Proof.)} Before starting, we need to observe that in a convex permutomino all the $\alpha$ and
$\gamma$ points belong to the permutation $\pi_1$, thus by  {\bf (z1)} they can also lie on the diagonal $x=y$;
on the contrary, the $\beta$ and $\delta$ points belong to $\pi_2$, then by  {\bf (z2)} all the $\beta$ (resp.
$\delta$) points must remain strictly above (resp. below) the diagonal $x+y=n+1$.

\medskip

\noindent $(\Longrightarrow)$ By Proposition \ref{und} we have that $\sigma$ is lower unimodal. Then, we have to
prove that $\pi$ may not be decomposed into the direct difference of two permutations, $\pi = \theta \ominus
\theta'$.

If $\pi(1)<\pi(n)$ the property is straightforward. Let us consider the case $\pi(1)>\pi(n)$, and assume that
$\pi = \theta \ominus \theta'$ for some permutations $\theta$ and $\theta'$. We will prove that if the vertices
of polygon $P$ define the permutation $\pi$, then the boundary of $P$ crosses itself, hence $P$ is not a
permutomino.

\smallskip

Let us assume that $P$ is a convex permutomino associated with $\pi = \theta \ominus \theta'$. We start by
observing that the $\beta$ and the $\delta$ paths of $P$ may not be empty. In fact, if the $\beta$ path is empty,
then $\pi (n)=n>\pi(1)$, against the hypothesis. Similarly, if the $\delta$ path is empty, then $\pi
(1)=1<\pi(n)$. Essentially for the same reason, both $\theta$ and $\theta'$ must have more than one element.

\begin{figure}[htb]
\begin{center}
\centerline{\hbox{\includegraphics[width=0.50\textwidth]{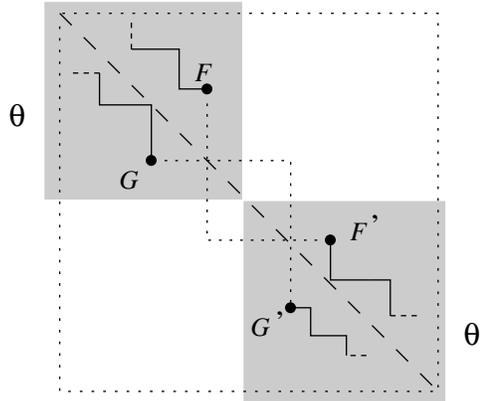}}} \caption{ If $\pi = \theta \ominus \theta'$
then the boundary of every polygon associated with $\pi$ crosses itself.}\label{th1}
\end{center}
\end{figure}

As we observed, the points of $\theta$ (resp. $\theta'$) in the $\beta$ path of $P$, are placed strictly above
the diagonal $x+y=n+1$. Let $F$ (resp. $F'$) be the rightmost (resp. leftmost) of these points. Similarly, there
must be at least one point of $\theta$ (resp. $\theta'$) in the $\delta$ path of $P$, placed strictly below the
diagonal $x+y=n+1$. Let $G$ (resp. $G'$) be the rightmost (resp. leftmost) of these points. The situation is
schematically sketched in Figure \ref{th1}.

Since $F$ and $F'$ are consecutive points in the $\beta$ path of $P$, they must be connected by means of a path
that goes down and then right, and, similarly, since $G'$ and $G$ are two consecutive points in the $\delta$
path, they must be connected by means of a path that goes up and then left. These two paths necessarily cross in
at least two points, and their intersections must be on the diagonal $x+y=n+1$.

\medskip

\begin{figure}[htb]
\begin{center}
\centerline{\hbox{\includegraphics[width=0.8\textwidth]{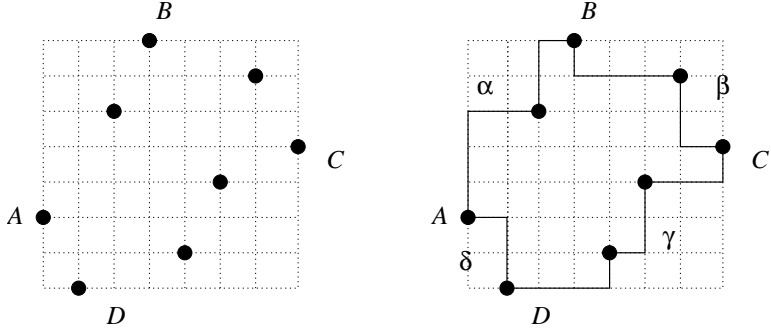}}} \caption{Given the permutation $\pi =
(3,1,6,8,2,4,7,5)$ satisfying conditions 1. and 2., we construct the $\alpha$, $\beta$, $\gamma$, and $\delta$
paths.}\label{th2}
\end{center}
\end{figure}

\noindent \noindent $(\Longleftarrow)$ Clearly condition 2. implies that $\pi (1) < n$ and $\pi (n)>1$, which are
necessary conditions for $\pi \in \widetilde{\mathcal C}_n$. We start building up a polygon $P$ such that
$\pi_1(P)=P$, and then prove that $P$ is a permutomino. As usual, let us consider the points
$$A=(1,\pi(1)), \quad B=(\pi^{-1}(n),n), \quad C=(n,\pi(n)),
\quad D=(\pi^{-1}(1),1).$$ The $\alpha$ path of $P$ goes from $A$ to $B$, and it is constructed connecting the
points of $\mu$ increasing sequence; more formally, if $\mu (i_l)$ and $\mu (i_{l+1})$ are two consecutive points
of $\mu$, with $\mu (i_l)<\mu (i_{l+1})\leq n$, we connect them by means of a path $$ 1^{\mu (i_{l+1}) - \mu
(i_l)} \, 0^{i_{l+1} -i_l} ,$$ (where $1$ denotes the vertical, and $0$ the horizontal unit step). Similarly we
construct the $\beta$ path, from $B$ to $C$, the $\gamma$ path from $C$ to $D$, and the $\delta$ path from $D$ to
$A$. Since the subsequence $\sigma$ is lower unimodal the obtained polygon is convex (see Figure \ref{th2}).

\begin{figure}[htb]
\begin{center}
\centerline{\hbox{\includegraphics[width=.8\textwidth]{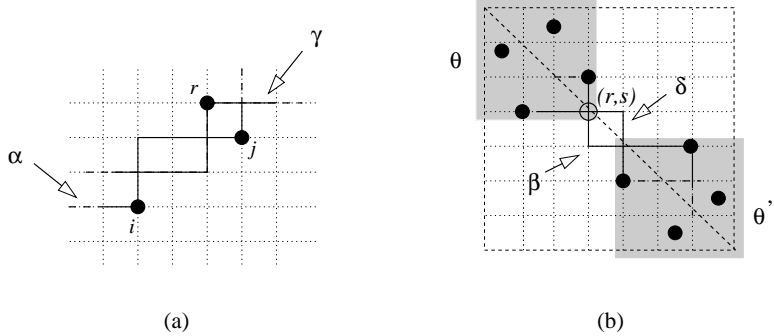}}} \caption{(a) The $\alpha$ path and the $\gamma$
path may not cross; (b) The $\beta$ path and the $\delta$ path may not cross.}\label{th3}
\end{center}
\end{figure}

Now we must prove that the four paths we have defined may not cross themselves. First we show that the $\alpha$
path and the $\gamma$ path may not cross. In fact, if this happened, there would be a point $(r,\pi(r))$ in the
path $\gamma$, and two points $(i,\pi(i))$ and $(j,\pi(j))$ in the path $\alpha$, such that $i<r<j$, and
$\pi(i)<\pi(r)>\pi(j)$ (see Figure \ref{th3} (a)). In this case, according to the definition, $\pi(r)$ should
belong to $\mu$, and then $(r,\pi(r))$ should be in the path $\alpha$, and not in $\gamma$.

Finally we prove that the paths $\beta$ and $\delta$ may not cross. In fact, if they cross, their intersection
should necessarily be on the diagonal $x+y=n+1$; if $(r,s)$ is the intersection point having minimum abscissa,
then the reader can easily check, by considering the various possibilities, that the points $(i,\pi(i))$ of $\pi$
satisfy:
$$i\leq r \; \; \mbox{ if and only if } \; \; \pi(i)\geq s $$
(see Figure \ref{th3} (b)). Therefore, setting
$$\theta =\left \{ (i,\pi(i)-s+1) : i\leq r \right \} $$
we have that $\theta$ is a permutation of ${\mathcal S}_r$, and letting
$$ \theta ' = \{ \, (i,\pi(i) \, : \, i>r \, \}$$
we see that $\pi  = \theta \ominus \theta '$, against the hypothesis. $\qed$

\smallskip

There is an interesting refinement of the previous general theorem, which applies to a particular subset of the
permutations of ${\mathcal S}_n$.

\begin{corollary}\label{caratt_conv1}
{\em Let $\pi \in {\mathcal S}_n$, such that $\pi(1)<\pi(n)$. Then $\pi \in \widetilde{\mathcal C}_n$ if and only
if $\sigma$ is lower unimodal.}
\end{corollary}

\medskip

\begin{figure}[htb]
\begin{center}
\centerline{\hbox{\includegraphics[width=.9\textwidth]{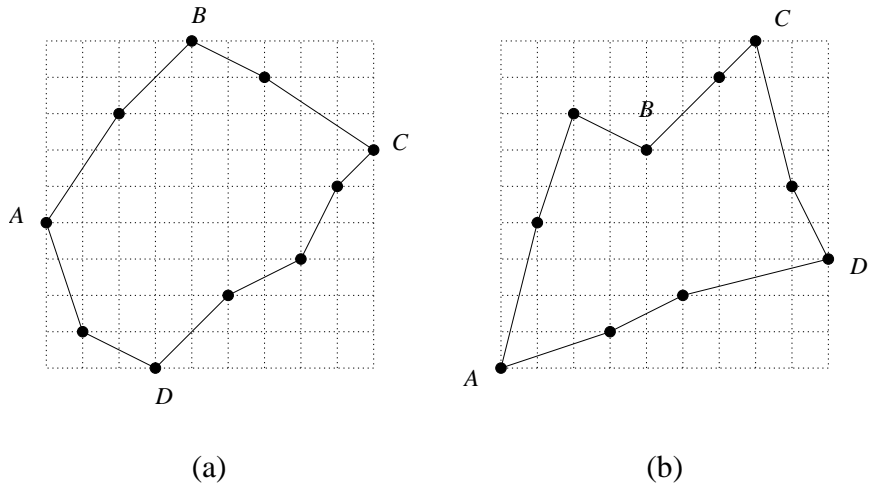}}} \caption{(a) a square permutation and the
associated $4$-face polygon; (b) a $4$ face polygon defined by a non square permutation.}\label{square}
\end{center}
\end{figure}
At the end of this section we would like to point out an interesting connection between the permutations
associated with convex permutominoes and another kind of combinatorial objects treated in some recent works. We
are referring to the so called {\em  $k$-faces permutation polygons} defined by T. Mansour and S. Severini in
\cite{toufik}. In order to construct a polygon from a given permutation $\pi$ in an unambiguous way, they find
the set of left-right minima and the set of right-left minima. An entry which is neither a left-right minimum nor
a right-left minimum is said to be a \emph{source}, together with the first and the last entry (which are also a
left-right minimum and a right-left minimum, respectively). Finally, two entries of $\pi$ are connected with an
edge if they are two consecutive left-right minima or right-left minima or sources. A maximal path of increasing
or decreasing edges defines a \emph{face}. If the obtained polygon has $k$ faces, than it is said to be a
$k$-$faces$ polygon. A permutation is said to be $square$ if the sequence of the sources lies in at most two
faces. The set of the square permutations of length $n$ is denoted by ${\mathcal Q}_n$. We note that a square
permutation has at most four faces, but the inverse statement does not hold: the permutation
$(1,5,8,2,7,3,9,10,6,4)$ has four faces and it is not square. Figure \ref{square} depicts an example.

Connecting all pairs of consecutive points of the sequences $\mu$ and $\sigma$ we obtain a polygon which may not
coincide with the polygon obtained from the definition of Mansour and Severini, as the reader can easily check
with the permutation $(1,2,4,3)$. It is however simple to state the following

\begin{proposition}\label{stopponi}
Given a permutation $\pi \in {\mathcal S}_n$, then $\pi \in {\mathcal Q}_n$ if and only if $\sigma(\pi)$ is lower
unimodal.
\end{proposition}

We point out that the square permutations coincide with the of the \emph{convex permutations}, introduced by
Waton \cite{waton}. In his PhD thesis the author characterizes the convex permutations in terms of forbidden
patterns. More precisely, he proves that the convex permutations are all the permutations avoiding the following
sixteen patterns of length five:
$$
\{52341, 52314, 51342, 51324, 42351, 42315, 41352, 41325$$ $$ 25341, 25314, 15342, 15324, 24351, 24315, 14352,
14325\}.
$$

\bigskip\noindent
All the relations between ${\mathcal Q}_n$, ${\mathcal C}_n$ and ${\mathcal C}'_n$ are exploited in the next
section, where, in particular, it is proved that, given a permutation $\pi$, then $\pi\in {\mathcal Q}_n$ if and
only if $\pi\in {\mathcal C}_n\cup{\mathcal C'}_n$.

Mansour and Severini \cite{toufik} (and independently Waton \cite{waton}) prove that the number $Q_{n+1}$ of
square permutation of size $n+1$ is
\begin{equation}
Q_{n+1}=2(n+3)4^{n-2}-4(2n-3){{2(n-2)}\choose {n-2}}, \label{squar}
\end{equation}
defining the sequence $1,2,6,24,104,464,2088, \ldots$ (A128652 in \cite{sloane}).

\subsection{The relation between the number of permutations and the number convex permutominoes}

Let $\pi \in \widetilde{\mathcal C}_n$, and $\mu$ and $\sigma$ defined as above. Let ${\mathcal F}(\pi)$ (briefly
${\mathcal F}$) denote the set of fixed points of $\pi$ lying in the increasing part of the sequence $\mu$ and
which are different from $1$ and $n$. We call the points in $\mathcal F$ the {\em free fixed points} of $\pi$.

For instance, concerning the permutation $\pi = (2,1,3,4,7,6,5)$ we have $\mu = ( 2,3,4,7,6,5)$, $\sigma =
(2,1,5)$, and ${\mathcal F}(\pi)=\{ 3,4 \}$; here $6$ is a fixed point of $\pi$ but it is not on the increasing
sequence of $\mu$, then it is not free. By definition, a permutation in $\widetilde{\mathcal C}_n$ can have no
free fixed points (e.g., the permutation associated with the permutomino in Figure~\ref{conv_new}), and at most
$n-2$ free fixed points (as the identity $(1, \ldots , n)$).

\begin{theorem}\label{car_2}
Let $\pi\in \widetilde{\mathcal C}_n$, and let ${\mathcal F}(\pi)$ be the set of free fixed points of $\pi$. Then
we have:
$$ \left | \; \left [ \; \pi \; \right ] \; \right | \; = \; 2^{\left | {\mathcal F}(\pi) \right |}.$$
\end{theorem}

\noindent {\bf (Proof.)} Since $\pi\in \widetilde{\mathcal C}_n$ there exists a permutomino $P$ associated with
$\pi$. If we look at the permutation matrix defined by the reentrant points of $P$, we see that all the free
fixed points of $\pi$ can be only of type $\alpha$ or $\gamma$, while the type of all the other reentrant points
of $\pi$ is established. It is easy to check that in any way we set the typology of the free fixed points in
$\alpha$ or $\gamma$ we obtain, starting from the matrix of $P$, a permutation matrix which defines a convex
permutomino associated with $\pi$, and in this way we get all the convex permutominoes associated with the
permutation $\pi$. $\qed$


\medskip

\begin{figure}[htb]
\begin{center}
\centerline{\hbox{\includegraphics[width=1\textwidth]{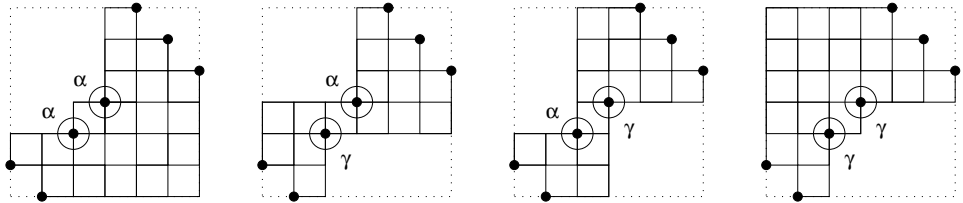}}} \caption{The four convex permutominoes
associated with the permutation $\pi = (2,1,3,4,7,6,5)$. The two free fixed points are encircled.}\label{classi}
\end{center}
\end{figure}

\noindent Applying Theorem \ref{car_2} we have that the number of convex permutominoes associated with $\pi =
(2,1,3,4,7,6,5)$  is $2^2=4$, as shown in Figure \ref{classi}. Moreover, Theorem \ref{car_2} leads to an
interesting property.

\begin{proposition}\label{caratt_conv2}
Let $\pi \in \widetilde{\mathcal C}_n$, with $\pi(1)>\pi(n)$. Then there is only one convex permutomino
associated with $\pi$, i.e., $\left | \, \left [ \, \pi \, \right ] \, \right | \, = \, 1$.
\end{proposition}

\smallskip

\noindent {\bf (Proof.)} If $\pi(1)>\pi(n)$ then all the points in the increasing part of $\mu$ are strictly
above the diagonal $x=y$, then $\pi$ cannot have free fixed points. The thesis is then straightforward. $\qed$

\medskip

Let us now introduce the sets $\widetilde{\mathcal C}_{n,k}$ of permutations having exactly $k$ free fixed
points, with $0\leq k \leq n-2$. We easily derive the following relations:
\begin{equation}\label{xx}
\widetilde{C}_{n} = \sum _{k=0}^{n-2} \; \left | \widetilde{\mathcal C}_{n,k} \right | \qquad \qquad \qquad
{C}_{n} = \sum _{k=0}^{n-2} \; 2^k \; \left | \widetilde{\mathcal C}_{n,k} \right |.
\end{equation}

\section{The cardinality of $\widetilde{\mathcal C}_{n}$}

In order to find a formula to express $\widetilde{C}_n$, it is now sufficient to count how many permutations of
${\mathcal Q}_n$ can be decomposed into the direct difference of other permutations. We say that a square
permutation is {\em indecomposable} if it is not the direct difference of two permutations. For any $k\geq 2$,
let $${\mathcal B}_{n,k}=\left \{ \, \pi \in {\mathcal Q}_n: \pi = \theta_1 \ominus \ldots \ominus \theta_k , \,
\theta_i \mbox{ indecomposable}, \, 1\leq i \leq k \, \right \}$$ be the set of square permutations which are
direct difference of exactly $k$ indecomposable permutations, and $$ {\mathcal B}_n = \bigcup _{k\geq 2}
{\mathcal B}_{n,k}.$$ For any $n,k\geq 2$, let ${\mathcal T}_{n,k}$ be the class of the sequences
$(P_1,\dots,P_k)$ such that:
\begin{description}
    \item[i] $P_1$ and $P_k$ are (possibly empty) directed convex permutominoes,
    \item[ii] $P_2, \dots, P_{k-1}$ are (possibly empty) parallelogram
    permutominoes,
\end{description}
 and such that the sum of the dimensions of $P_1,\dots,P_k$ is equal to
 $n$.

\begin{proposition}\label{filippo}
There is a bijective correspondence between the elements of ${\mathcal B}_{n,k}$ and the elements of ${\mathcal
T}_{n,k}$, so that the two classes have the same cardinality.
\end{proposition}

\noindent {\bf (Proof.)} Let us consider $(P_1,\dots,P_k)\in {\mathcal T}_{n,k}$, we construct the corresponding
permutation $\pi=\delta_1 \ominus \dots \ominus \delta_k$ as follows. For any $1\leq i \leq k$, if $P_i$ is the
empty permutomino, then $\delta_i=(1)$, otherwise:

\begin{description}
    \item[i] for all $i$ with $1\leq i\leq k-1$, $\delta_i$ is the reversal
    of $\pi_2(P_i)$ (i.e., the permutation $\pi_1$ associated with the symmetric
    permutomino of $P_i$ with respect to the $y$- axis).
    \item[ii] $\delta_k$ is the complement of $\pi_2(P_k)$ (i.e., it is the
    permutation $\pi_1$ associated with the symmetric permutomino of $P_k$ with
    respect to the $x$-axis).
\end{description}

\begin{figure}[htb]
\begin{center}
\centerline{\hbox{\includegraphics[width=1.1\textwidth]{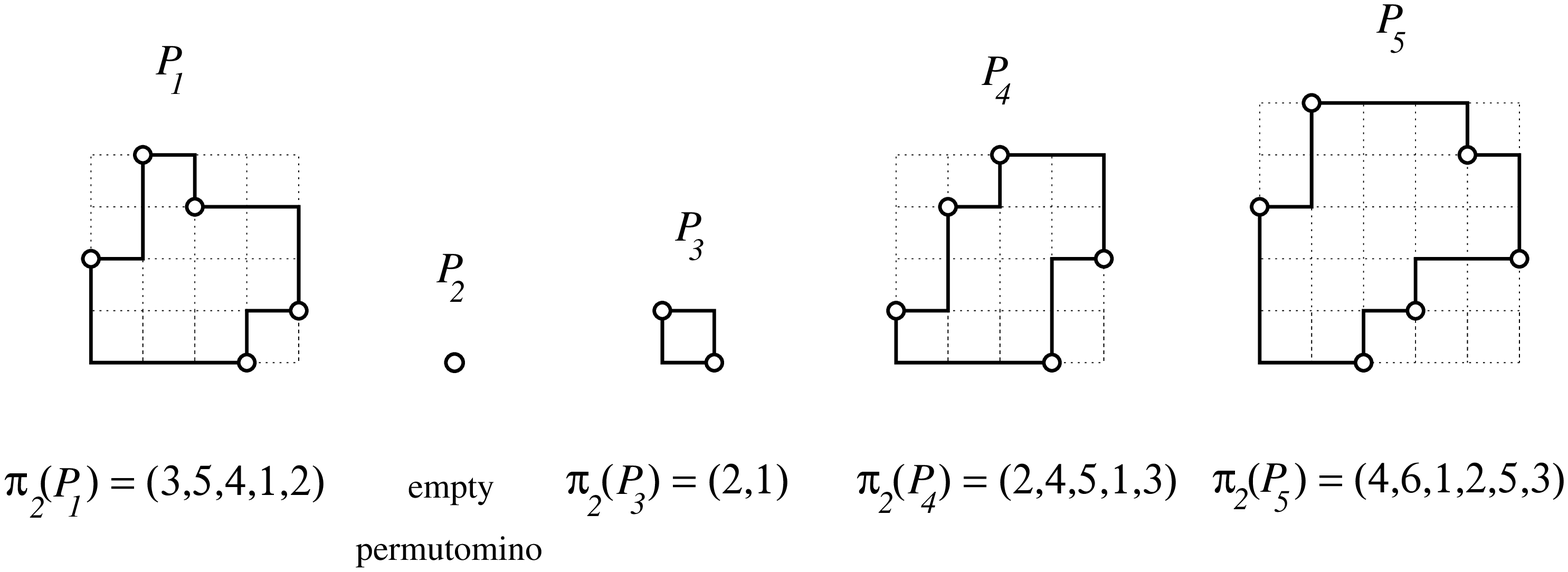}}} \caption{An element of
 ${\mathcal T}_{19,5}$, constituted of a sequence of five permutominoes,
 and the associated permutations.}\label{squaredd}
\end{center}
\end{figure}

\begin{figure}[htb]
\begin{center}
\centerline{\hbox{\includegraphics[width=1.1\textwidth]{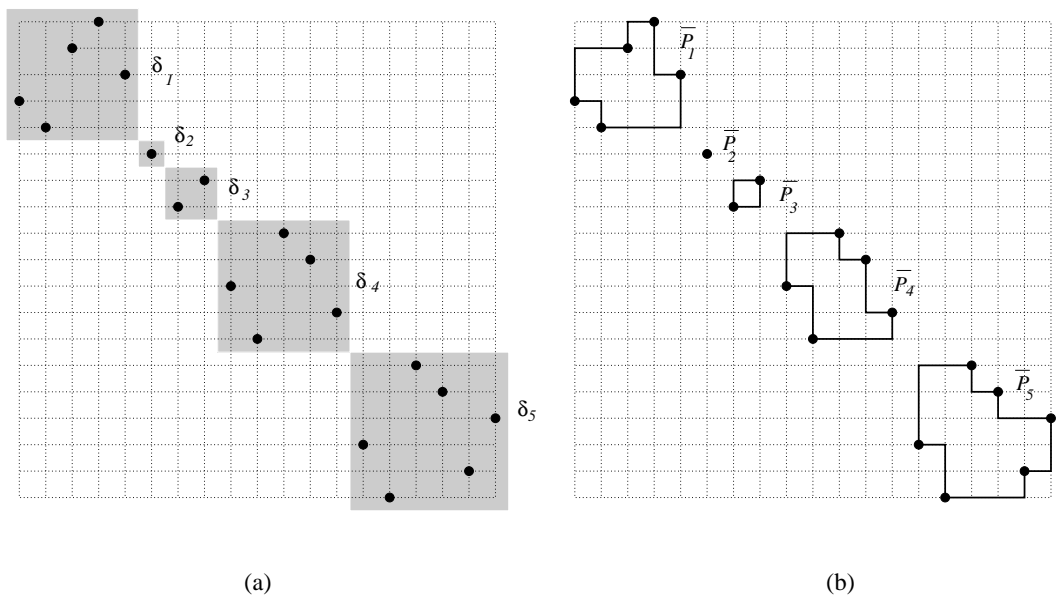}}} \caption{(a) a square permutation which
can be decomposed into the direct difference of five indecomposable permutations; (b) the five permutominoes
associated with them. For each permutomino $P_i$, we denote by $\bar{P}_i$ the corresponding reflected
permutomino.} \label{squared}
\end{center}
\end{figure}

For example, starting from the sequence of permutominoes in Figure \ref{squaredd} we obtain the following
permutations: $\delta_1=(2,1,4,5,3)$ is obtained from the permutomino $P_1$ such that $\pi_2(P_1)=(3,5,4,1,2)$;
$\delta_2=(1)$ is obtained from the empty permutomino $P_2$; $\delta_3=(1,2)$ is obtained from $P_3$;
$\delta_4=(3,1,5,4,2)$ is obtained from the permutomino $P4$ such that $\pi_2(P_4)=(2,4,5,1,3)$. Moreover,
$\delta_5=(3,1,6,5,2,4)$ is the complement of $\tau=(4,6,1,2,5,3)$ which is such that $\pi_2(P_5)=\tau$. Then, as
showed in Figure \ref{squared} we obtain the permutation
$\pi=\delta_1\ominus\delta_2\ominus\delta_4\ominus\delta_4\ominus\delta_5$,
$$\pi=(16,15,18,19,17,14,12,13,9,7,11,10,8,3,1,6,5,2,4)$$
We note that the points in the increasing part of $\mu(\pi)$ are precisely the points of the increasing part of
$\mu(\delta_1)$; the points in the increasing part of $\sigma(\pi)$ are the points of the increasing part of
$\sigma(\delta_k)$; the points in the decreasing part of $\mu(\pi)$ are given by the sequence of points of the
decreasing parts of $\mu(\delta_1),\dots,\mu(\delta_k)$; finally, the points in the decreasing part of
$\sigma(\pi)$ are given by the sequence of the points of the decreasing parts of
$\sigma(\delta_1),\dots,\sigma(\delta_k)$. Then, we have that $\pi\in\mathcal Q_n$ and then $\pi\in {\mathcal
B}_{n,k}$.

\medskip

Conversely, let $\pi\in {\mathcal B}_{n,k}$, with $\pi=\delta_1 \ominus \dots \ominus \delta_k$. By the previous
considerations we have that $\pi\in {\mathcal Q}_n$, and then it is clear that, for each component $\delta _i$,
the sequence $\mu(\delta_i)$ is upper unimodal, and $\sigma(\delta_i)$ is lower unimodal.

If $\delta _i$ is the one element permutation, then it is associated with the empty permutomino. Otherwise, if a
permutation $\delta_i$ is indecomposable and has dimension greater than $1$ it is clearly associated with a
polygon with exactly one side for every abscissa and ordinate and with the border which does not intersect
itself. These two conditions are sufficient to state that $\delta_i$ is associated with a convex permutomino, and
in particular the reader can easily observe the following properties, due to its the indecomposability:

\begin{enumerate}
    \item there is exactly one directed convex permutomino $P_1$
corresponding to $\delta_1$, and it is the reflection according to the $y$-axis of a permutomino associated with
$\delta _1$;
    \item for any $2\leq i\leq k-1$, there is exactly one parallelogram
permutomino $P_i$ corresponding to $\delta_i$, and it is the reflection according to the $y$-axis of a
permutomino associated with $\delta _i$;
    \item there is exactly one directed convex permutomino $P_k$
corresponding to $\delta_k$, and it is the reflection according to the $x$-axis of a permutomino associated with
$\delta _k$.
\end{enumerate}
We have thus the sequence $(P_1,\dots,P_k)\in {\mathcal T}_{n,k}$. $\qed$

\bigskip

If we denote by $B_n$ (resp. $B_{n,k}$) the cardinality of ${\mathcal B}_n$ (resp. ${\mathcal B}_{n,k}$), by
Proposition \ref{stopponi} we have
$$ \widetilde{C}_n = Q_n - B_n.$$

\noindent Let us pass to generating functions, denoting by:

\begin{enumerate}
    \item $P(x)$ (resp. $D(x)$) the generating function of
    parallelogram permutominoes (resp. $P(x)$), hence
    \begin{eqnarray*}
     P(x) = \frac{1-\sqrt{1-4x}}{2} &=& x+x^2+2x^3+5x^4+14x^5+
    \ldots  \\
     D(x) = \frac{x}{2} \left( \frac{1}{\sqrt{1-4x}} +1 \right ) &=& x+x^2+3x^3+10x^4+35x^5+
    \ldots ;
    \end{eqnarray*}
    \item $B_k(x)$ (resp. $B(x)$) the generating function of the
    numbers $\{ B_{k,n} \}_{n\geq 0}$, $k\geq 2$ (resp. $\{ B_{n} \}_{n\geq
    0}$).
\end{enumerate}

\noindent Due to Proposition  \ref{filippo}, for any $k\geq 2$, we have that $B_k(x)=D^2(x)P^{k-2}(x)$ and then
$$ B(x)=\sum_{k\geq 0} D^2(x) P^{k-2}(x)
=\frac{D^2(x)}{1-P(x)}=\frac{1}{2}\left ( \frac{x^2}{1-4x} + \frac{x^2}{\sqrt{1-4x}} \right ).$$ Therefore
$$B_{n+2}=\frac{1}{2}\left ( 4^{n} +
{{2n}\choose n} \right )=\sum_{i=0}^n {{2n}\choose i}. $$ Now it is easy to determine the cardinality of
$\widetilde{\mathcal C}_{n}$. For simplicity of notation we will express most of the following formulas in terms
of $n+1$ instead of~$n$.
\begin{proposition}\label{wwww}
{\em The number of permutations of $\widetilde{\mathcal C}_{n+1}$ is
\begin{equation}\label{solution}
 \, 2 \, (n+2) \, 4^{n-2} \, - \, \frac{n}{4} \, \left ( {\frac{3-4n}{1-2n}}
\right ) \, {{2n} \choose {n}}, \qquad n \geq 1.
\end{equation}}
\end{proposition}

\noindent {\bf (Proof.)} In fact, for any $n\geq 2$, we have $\widetilde{C}_{n}=Q_{n}-B_{n}$, then the result is
straightforward. $\qed$

\bigskip\noindent
For the sake of completeness, in Table \ref{tab} we list the first terms of the sequences involved in the
preceding formulas.

\begin{table}[htb]
\begin{center}
\begin{tabular}{l|ccccccccc}
  sequence &$1$ &$2$ &$3$ &$4$ &$5$ &$6$ &$7$ &$8$ &$\ldots$ \\
  \hline
    \\
  $Q_n$ &$1$ &$2$ &$6$ &$24$ &$104$ &$464$ &$2088$ &$9392$ &$\ldots$ \\
    \\
  $B_n$ & &$1$ &$3$ &$11$ &$42$ &$163$ &$638$ &$2510$ &$\ldots$ \\
    \\
  $\widetilde{C}_n $ &$1$ &$1$ &$3$ &$13$ &$62$ &$301$ &$1450$ &$6882$
&$\ldots$ \\
\end{tabular}
\caption{The first terms of the sequences $Q_n$, $B_n$, $\widetilde{C}_n $, starting with $n=1$.}\label{tab}
\end{center}
\end{table}

\medskip

In ending the paper we would like to point out some other results that directly come out from the one stated in
Proposition \ref{wwww}. First we observe that the number of permutations $\pi \in \widetilde{\mathcal C}_{n}$ for
which $\pi (1) < \pi(n)$ is equal to $\frac{1}{2} Q_n$, while the number of those for which $\pi (1)
> \pi(n)$ is equal to
$$ \frac{1}{2}Q_n - B_n = \widetilde{C}_{n} - \frac{1}{2}Q_n ,$$
and the $(n+1)$th term of this difference is equal to
\begin{equation}\label{exp}
(n+1)4^{n-2}-\frac{n}{2}{{2n+1}\choose {n-1}},
\end{equation}
whose first terms are $1,10,69,406, 2186, 11124,\ldots$, (sequence A038806 in \cite{sloane}).

\medskip

\noindent Moreover,
it is also possible to consider the set $\widetilde{\mathcal C}_{n} \cap \widetilde{\mathcal C}'_{n}$, i.e., the
set of the permutations $\pi$ for which there is at least one convex permutomino $P$ such that $\pi_1(P)=\pi$ and
one convex permutomino $P'$ such that $\pi_2(P')=\pi$. For instance, we have:
$$
\begin{array}{ll}
\widetilde{\mathcal C}_{3} \cap \widetilde{\mathcal C}'_{3}=\emptyset, \\
\\
\widetilde{\mathcal C}_{4} \cap \widetilde{\mathcal C}'_{4}=\{ (2,4,1,3), (3,1,4,2) \}.
\end{array}
$$
We start by recalling that $\pi \in \widetilde{\mathcal C}_{n}$ if and only if $\pi^R \in \widetilde{\mathcal
C}'_{n}$
.

\begin{proposition}\label{opo}
{\em A permutation $\pi \in {\mathcal Q}_n$ if and only if $\pi \in \widetilde{\mathcal C}_{n} \cup
\widetilde{\mathcal C}'_{n}$.}
\end{proposition}

\noindent {\bf (Proof.)} ($\Leftarrow$) If $\pi$ is a square permutation but it is not in $\widetilde{\mathcal
C}_{n}$, then necessarily $\pi (1) > \pi (n)$. Hence, if we consider $\pi^M$, we have $\pi^M(1)<\pi^M(n)$, and
$\pi^M\in \widetilde{\mathcal C}_{n}$, then $\pi\in \widetilde{\mathcal C}'_{n}$.

\smallskip

\noindent ($\Rightarrow$) Trivial. $\qed$

\medskip

Finally, since $\left | \widetilde{\mathcal C}'_n \right |= \left | \widetilde{\mathcal C}_n \right |$, and $
Q_{n} = 2\widetilde{C}_n - \left | \, \widetilde{{\mathcal C}}_n \cap \widetilde{{{\mathcal C}'}}_n \right |$, we
can state the following.

\begin{proposition}\label{opoz}
{\em For any $n \geq 2$, we have
\begin{equation}\label{ds} \left |
\, \, \widetilde{{\mathcal C}}_{n} \cap \widetilde{{{\mathcal C}'}}_{n} \right | = \widetilde{{C}}_{n} -
B_{n}=Q_n-2B_n.
\end{equation}}
\end{proposition}

The reader can easily recognize that the numbers defined by (\ref{ds}) are the double of the ones expressed by
the formula in (\ref{exp}), so that

\begin{equation}\label{gek}
\left | \, \widetilde{{\mathcal C}}_{n+1} \cap \widetilde{{{\mathcal C}'}}_{n+1} \right | =
2(n+1)4^{n-2}-{{2n-1}\choose {n-1}}.
\end{equation}

\section{Further work}
Here we outline the main open problems and research lines on the class of permutominoes.
\begin{enumerate}
    \item It would be natural to look for a combinatorial proof of the formula (\ref{co}) for the number of
    convex permutominoes
    and (\ref{solution}) for the number of permutations associated with convex permutominoes. These proofs could be
    obtained using the matrix characterization for convex permutominoes provided in Section \ref{rif}.

    \item The main results of the paper have been obtained in an analytical
    way. In particular from (\ref{co}) and (\ref{solution}) we
    have a direct relation between convex permutominoes and
    permutations, obtaining
    \begin{equation}\label{inter}
    C_{n+2} = \widetilde{C}_{n+2} + \frac{1}{2} \left ( 4^n - {{2n} \choose n} \right
    ),
    \end{equation}
    which requires a combinatorial explanation. In particular,
    recalling that
    $$ C_n = \sum _{\pi \in \widetilde{\mathcal C}_{n}} \left | [  \pi ]\right
    |, $$
    the right term of (\ref{inter}) is the number of convex
    permutominoes which are determined by the permutations having
    at least one free fixed point.

    Moreover, from (\ref{wwww}) and (\ref{inter}) we get that
    $$ Q_{n+2} = C_{n+2} + {{2n} \choose n}, $$
    and also this identity cannot be clearly explained using the
    combinatorial arguments used in the paper.

    \medskip

    From (\ref{gek}) we have that the generating function of the
    permutations in $\widetilde{{\mathcal C}}_{n} \cap \widetilde{{{\mathcal
    C}'}}_{n}$ is
    $$ 2\left ( \frac{x^2 c(x)}{1-4x} \right ) ^2, $$
    where $c(x)$ denotes the generating function of Catalan
    numbers. While the factor $2$ can be easily explained,
    since for any $\pi \in \widetilde{{\mathcal C}}_{n} \cap \widetilde{{{\mathcal
    C}'}}_{n}$, also $\pi^M \in \widetilde{{\mathcal C}}_{n} \cap \widetilde{{{\mathcal
    C}'}}_{n}$, and clearly $\pi \neq \pi '$, the convolution of
    Catalan numbers and the powers of four begs for a
    combinatorial interpretation.

    \item
    We would like to consider the characterization and the enumeration of the permutations associated with other
    classes of permutominoes, possibly including the class of convex permutominoes. For instance, if we take the
    class of column convex permutominoes, we observe that Proposition \ref{spigoli} does not hold.
    In particular, one can see that, if the permutomino is not convex,
    then the set of reentrant points does not form a permutation matrix (Figure \ref{ccc}).

    \begin{figure}[htb]
    \begin{center}
    \centerline{\hbox{\includegraphics[width=1\textwidth]{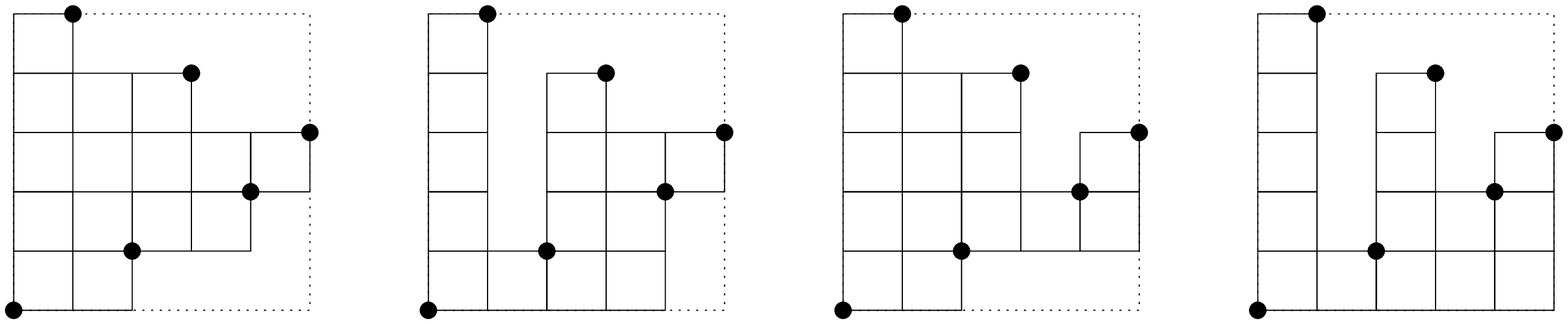}}}
    \caption{The four column convex permutominoes associated with the permutation $(1,6,2,5,3,4)$; only the
    leftmost is convex}\label{ccc}
    \end{center}
    \end{figure}

    Moreover, it might be interesting to determine an extension of Theorem \ref{car_2} for the class of
    column convex permutominoes, i.e., to characterize the set of column convex permutominoes
    associated with a given permutation. For instance, we observe that while there is one convex
    permutomino associated with $\pi=(1,6,2,5,3,4)$, there are four column convex permutominoes
    associated with $\pi$ (Figure \ref{ccc}).

\end{enumerate}


\begin{thebibliography}{20}


\bibitem{bmrr} G. Barequet, M. Moffie, A. Ribó, G. Rote, Counting Polyominoes on Twisted Cylinders,
Proc. of  {\em EuroComb '05, European Conference on Combinatorics, Graph Theory and Applications, Ed. S. Felsner,
Disc. Math. Theor. Comput. Sci.} Proceedings AE,  369-374.

\bibitem{Gi}
D. Beauquier, M. Nivat, Tiling the plane with one tile, Proc. of {\em 6th Annual Symposium on Computational
geometry}, Berkeley, CA, ACM press (1990) 128-138.


\bibitem{milanesi}
Boldi, P., Lonati, V., Radicioni, R., Santini, M.: {\rm The number of convex permutominoes.}, Proc. of {\it LATA
2007, International Conference on Language and Automata Theory and Applications}, Tarragona, Spain, (2007).

\bibitem{mbm2}
M. Bousquet-M\`elou, A method for the enumeration of various classes of column convex polygons, {\it Disc. Math.}
154 (1996) 1--25.

\bibitem{mbm}
M. Bousquet-M\`elou, A. J. Guttmann, Enumeration of three dimensional convex polygons, {\it Ann. of Comb.} 1
(1997) 27--53.

\bibitem{gutman}
Brak, R., Guttmann, A. J., Enting, I. G.: Exact solution of the row-convex polygon perimeter generating function,
{\it J. Phys.} A 23 (1990) L2319--L2326.

\bibitem{brlek}
Brlek, S., Labelle, G., Lacasse, A.: A Note on a Result of Daurat and Nivat, {\it Lecture Notes in Computer
Science}, Springer Berlin/Heidelberg, Vol. 3572 (2005) 189-198.

\bibitem{chang}
Chang, S.J., Lin, K.Y.: Rigorous results for the number of convex polygons on the square and honeycomb lattices,
{\em J. Phys. A: Math. Gen.} {\bf 21} (1988) 2635-2642.

\bibitem{mbintr1}
Conway, J.H., Lagarias, J.C.: Tiling with polyominoes and combinatorial group theory, {\it J. Comb. Th. A} {\bf
53} (1990) 183-208.

\bibitem{daurat}
Daurat, A., Nivat, M.: Salient and reentrant points of discrete sets, {\em Disc. Appl. Math.} {\bf 151} (2005)
106-121.

\bibitem{DV}
Delest, M., Viennot, X.G.: Algebraic languages and polyominoes enumeration, {\em Theor. Comp. Sci.} {\bf 34}
(1984) 169-206.

\bibitem{DDFR}
Del~Lungo, A., Duchi, E., Frosini, A., Rinaldi, S.: On the generation and enumeration of some classes of convex
polyominoes, {\em El. J. Comb.}, {\bf 11} (2004), \#R60.

%

\bibitem{rinaldi}
Disanto, F., Frosini, A., Pinzani, R., Rinaldi, S.: A closed formula for the number of convex permutominoes, {\em
ArXiv Mathematics e-prints} {math/0702550} (2007).


\bibitem{fanti}
Fanti, I., Frosini, A., Grazzini, E., Pinzani, R., Rinaldi, S.:
\newblock {\rm Polyominoes determined by permutations},
\newblock {\it (submitted)}.

\bibitem{Ga}
M. Gardner, Mathematical games,  {\em Scientific American}, (1958) Sept. 182--192, Nov. 136-142.

\bibitem{Go}
S. W. Golomb, {\em Polyominoes: Puzzles, Patterns, Problems, and Packings}, Princeton Academic Press, 1996.

\bibitem{dbintr39}
S. W. Golomb, Checker boards and polyominoes, {\em Amer. Math. Monthly}, {\bf 61} (1954) 675-682.

\bibitem{jensen-guttmann}
I. Jensen, A. J. Guttmann, Statistics of lattice animals (polyominoes) and polygons, {\em J. Phys. A}, {\bf 33}
(2000) 257-263.

\bibitem{incitti}
Incitti, F., Permutation diagrams, fixed points and Kazdhan-Lusztig $R$-polynomials, {\it Ann. Comb.}, {\bf 10},
N.3, (2006)  369-387.

\bibitem{toufik}
Mansour, T., Severini, S., Grid polygons from permutations and their enumeration by the kernel method, {\em ArXiv
Mathematics e-prints} {math/0603225} (2006).

\bibitem{sloane}
Sloane, N.J.A.; The On-Line Encyclopedia of Integer Sequences, {\it http://www.research.att.com/ $\sim$
njas/sequences/}

\bibitem{stan}
R. P. Stanley,
\newblock {\it Enumerative Combinatorics,} Vol. 2,
\newblock {\rm Cambridge University Press,} Cambridge (1999).

\bibitem{waton}
Waton, S., On Permutation Classes Generated by Token Passing Networks, Gridding Matrices and Pictures: Three
Flavours of Involvement, {\em PhD Thesis}, University of St Andrews (2007).


\end{thebibliography}
\end{document}